\documentclass[a4paper,12pt]{article}
\usepackage{comment}
\usepackage{cite}
\usepackage{amsmath,amssymb,amsfonts}
\usepackage[T1]{fontenc}
\usepackage[utf8]{inputenc}
\usepackage{graphicx}
\usepackage{fancyhdr}
\usepackage{float}
\usepackage{xcolor}
\usepackage{authblk}
\usepackage{mathrsfs}
\usepackage{empheq}
\usepackage[hyphens]{url}
\usepackage{hyperref}
\usepackage{amsthm}
\usepackage{tikz}
\usetikzlibrary{decorations.markings}
\usepackage{enumitem}
\usepackage{geometry}

\theoremstyle{plain}

\theoremstyle{definition}

\theoremstyle{remark}

\begin{document}

\title{Elliptic integral identities derived from Coxeter's integrals}

\author[$\dagger$]{Jean-Christophe {\sc Pain}\footnote{jean-christophe.pain@cea.fr}\\
\small
$^1$CEA, DAM, DIF, F-91297 Arpajon, France\\
$^2$Universit\'e Paris-Saclay, CEA, Laboratoire Mati\'ere en Conditions Extrêmes,\\ 
F-91680 Bruy\`eres-le-Ch\^atel, France
}

\date{}

\maketitle

\begin{abstract}
We revisit the classical integrals introduced by Coxeter, not to recalculate their well-known exact values, but to use them as a tool to derive elliptic integral identities. By embedding Coxeter's first integral into a one-parameter family
$$
I(\lambda)=\int_{0}^{\pi/2}
\arccos\!\left(\frac{\cos\theta}{1+\lambda\cos\theta}\right)\,d\theta,
$$
and differentiating with respect to the parameter \(\lambda\), we show that the derivative $I'(\lambda)$ can be expressed as an elliptic-type integral. Integrating $I'(\lambda)$ between 0 and 2 yields the identity
$$
\int_0^2 \int_0^{\pi/2}
\frac{\cos^2\theta}
{(1+s\cos\theta)\sqrt{(1+s\cos\theta)^2-\cos^2\theta}}
\,d\theta\, ds=A-B=\frac{\pi^2}{12},
$$
where $A$ and $B$ are the first two so-called Coxeter integrals
$$
A = \int_0^{\pi/2} \arccos\!\left(\frac{\cos\theta}{1+2\cos\theta}\right) d\theta,
$$
and
$$
B = \int_0^{\pi/2} \arccos\!\left(\frac{1}{1+2\cos\theta}\right) d\theta.
$$
The derivative $I'(\lambda)$ can be expressed in terms of incomplete elliptic integrals of the first kind $F$ and of the third kind $\Pi$. This approach establishes a direct connection between classical Coxeter integrals and elliptic functions. The method highlights how well-known trigonometric integrals can serve as a bridge to explore properties and relations of elliptic integrals, offering new analytic insights beyond the original Coxeter evaluations.
\end{abstract}

\section{Introduction}

In a short note that has since become classical, Coxeter studied the following three remarkable integrals \cite{Coxeter1937,Vidiani2003}:
\[
A = \int_0^{\pi/2} \arccos\!\left(\frac{\cos\theta}{1+2\cos\theta}\right) d\theta,
\]
\[
B = \int_0^{\pi/2} \arccos\!\left(\frac{1}{1+2\cos\theta}\right) d\theta,
\]
and
\[
C = \int_0^{\pi/2} \arccos\!\left(\frac{1-\cos\theta}{2\cos\theta}\right) d\theta.
\]
These integrals have the striking property of taking values that are rational multiples of $\pi^2$:
\[
A = \frac{5\pi^2}{24}, 
\qquad 
B = \frac{\pi^2}{8}, 
\qquad 
C = \frac{11\pi^2}{72}.
\]
Beyond their analytic appearance, these formulas admit a deep geometric interpretation. They are closely related to the geometry of spherical and hyperbolic tetrahedra and to the differential formula of Schl\"afli \cite{Schlafli1858,Vidiani2003}, which relates the variation of the volume of a polyhedron (in geometries of constant curvature) to the variation of its dihedral angles:
\[
dV = -\frac{1}{2} \sum_{\text{edges}} \ell_i \, d\alpha_i.
\]
In this framework, Coxeter's integrals arise naturally as integrals of dihedral angles along specific deformations, and their reduction to rational multiples of $\pi^2$ reflects the highly symmetric structure of the underlying geometric configurations. Thus, the integrals $A$, $B$, and $C$ form a remarkable meeting point between elementary analysis, spherical geometry, and the theory of polyhedra, elegantly illustrating how seemingly innocent trigonometric expressions may conceal a deep geometric structure governed by Schläfli's formula \cite{Coxeter1968,Coxeter1973}.

Different methods were proposed for the calculation of such integrals. One of them relates to the so-called Ahmed integral \cite{Ahmed2002,Borwein2004,Nahin2015}. In this work, we present a novel perspective on the classical Coxeter integrals, not by recomputing their known exact values, but by using them as a tool to generate new relations between trigonometric integrals and elliptic integrals. We embed one of these integrals into the one-parameter family
\[
I(\lambda) = \int_0^{\pi/2} \arccos\!\left(\frac{\cos\theta}{1+\lambda \cos\theta}\right)\, d\theta,
\]
and show that differentiating with respect to the parameter \(\lambda\) leads naturally to an expression of the elliptic type. Integration then produces an identity in terms of an integral of $I'(\lambda)$ over the interval \([0,2]\), which is equal to the difference $A-B$ between the first two Coxeter integrals.  

The paper is organized as follows. In Section~2, we introduce the one-parameter family $I(\lambda)$ and compute its derivative with respect to the parameter. After a careful reduction based on the Weierstrass substitution, we obtain an explicit quartic representation of $I'(\lambda)$, exhibiting its elliptic structure. In Section~3, we integrate $I'(\lambda)$ over the interval $[0,2]$, provide an expression in terms of incomplete elliptic integrals $F$ and $\Pi$, and show how the classical evaluations of Coxeter's first two integrals (\textit{stricto sensu}, especially the first one $A$ actually) enter naturally as boundary values of the parametric family, leading to the identity
\[
\int_0^2 \int_0^{\pi/2}
\frac{\cos^2\theta}
{(1+s\cos\theta)\sqrt{(1+s\cos\theta)^2-\cos^2\theta}}
\, d\theta \, ds
=
\frac{\pi^2}{12}.
\]
We also briefly comment on the third Coxeter integral and discuss how this framework may be extended to other parametric deformations.

\section{A parametric family of integrals}

We consider the family
\[
I(\lambda)=\int_{0}^{\pi/2} \arccos\left(\frac{\cos\theta}{1+\lambda\cos\theta}\right)d\theta, \quad \lambda > -1.
\]
Let us set
\[
u(\theta,\lambda)=\frac{\cos\theta}{1+\lambda\cos\theta}.
\]
Then, the derivative is given by
\[
\frac{\partial}{\partial \lambda}I(\lambda)=\int_{0}^{\pi/2} \frac{-1}{\sqrt{1-u^2}} \frac{\partial u}{\partial \lambda} d\theta,
\]
where
\[
\frac{\partial u}{\partial \lambda}=-\frac{\cos^2\theta}{(1+\lambda\cos\theta)^2}.
\]
We thus have
\[
I'(\lambda)=\int_{0}^{\pi/2} \frac{\cos^2\theta}{(1+\lambda\cos\theta)\sqrt{(1+\lambda\cos\theta)^2-\cos^2\theta}} d\theta,
\]
as well as the identity
\[
(1+\lambda\cos\theta)^2-\cos^2\theta = 1+2\lambda\cos\theta+(\lambda^2-1)\cos^2\theta.
\]
Let us set
\[
a=\lambda^2-1.
\]
Using the Weierstrass substitution
\[
t=\tan\frac{\theta}{2}, \quad \cos\theta=\frac{1-t^2}{1+t^2}, \quad d\theta=\frac{2dt}{1+t^2},
\]
we obtain
\[
\cos^2\theta = \frac{(1-t^2)^2}{(1+t^2)^2}.
\]
Let us define the expression under the square root in the denominator as $\Delta(\theta)$. Substituting $\cos\theta$, we have
\[
\Delta(\theta)=1+2\lambda\frac{1-t^2}{1+t^2}+a\frac{(1-t^2)^2}{(1+t^2)^2}.
\]
Bringing everything over the common denominator $(1+t^2)^2$, the numerator $Q(t)$ becomes
\[
Q(t)=(1+t^2)^2+2\lambda(1-t^2)(1+t^2)+a(1-t^2)^2.
\]
We expand carefully:
\[
(1+t^2)^2 = 1+2t^2+t^4,
\]
\[
2\lambda(1-t^2)(1+t^2) = 2\lambda-2\lambda t^4,
\]
\[
a(1-t^2)^2 = a(1-2t^2+t^4).
\]
Collecting powers of $t$ gives, for the constant term:
\[
1+2\lambda+a = 1+2\lambda+(\lambda^2-1) = \lambda^2+2\lambda,
\]
for the coefficient of $t^2$:
\[
2-2a = 2-2(\lambda^2-1) = 4-2\lambda^2,
\]
and for the coefficient of $t^4$:
\[
1-2\lambda+a = 1-2\lambda+\lambda^2-1 = \lambda^2-2\lambda.
\]
We thus obtain
\[
\Delta(\theta) = \frac{(\lambda^2+2\lambda)+(4-2\lambda^2)t^2+(\lambda^2-2\lambda)t^4}{(1+t^2)^2},
\]
and therefore
\[
\sqrt{\Delta(\theta)} = \frac{\sqrt{(\lambda^2+2\lambda)+(4-2\lambda^2)t^2+(\lambda^2-2\lambda)t^4}}{1+t^2}.
\]
The integral becomes
\[
I'(\lambda) = \int_{0}^{1} \frac{\frac{(1-t^2)^2}{(1+t^2)^2} \cdot \frac{2}{1+t^2}}{\left(\frac{(1+t^2)+\lambda(1-t^2)}{1+t^2}\right) \frac{\sqrt{Q(t)}}{1+t^2}} dt,
\]
where $Q(t)$ is the quadratic polynomial in $u=t^2$. Cancelling the factors of $(1+t^2)$ yields
\[
I'(\lambda) = 2 \int_{0}^{1} \frac{(1-t^2)^2}{(1+t^2)\left((1+\lambda)+(1-\lambda)t^2\right)\sqrt{Q(t)}} dt.
\]
The discriminant of the polynomial $Q(u)=(\lambda^2-2\lambda)u^2+(4-2\lambda^2)u+(\lambda^2+2\lambda)$ is
\[
\Delta_u = (4-2\lambda^2)^2-4(\lambda^2-2\lambda)(\lambda^2+2\lambda).
\]
Since $(\lambda^2-2\lambda)(\lambda^2+2\lambda)=\lambda^2(\lambda-2)(\lambda+2)=\lambda^2(\lambda^2-4)$, we obtain
\[
\Delta_u = (16-16\lambda^2+4\lambda^4)-4(\lambda^4-4\lambda^2) = 16.
\]
The roots $u_{\pm}$ are therefore real and remarkably simple:
\[
u_{\pm} = \frac{2\lambda^2-4\pm 4}{2(\lambda^2-2\lambda)}.
\]
We thus obtain, for $\lambda\ne 0,1$, the exact representation
\[
I'(\lambda) = 2 \int_{0}^{1} \frac{(1-t^2)^2}{(1+t^2)\left((1+\lambda)+(1-\lambda)t^2\right)\sqrt{(\lambda^2+2\lambda)+(4-2\lambda^2)t^2+(\lambda^2-2\lambda)t^4}} dt,
\]
which is the correct elliptic quartic representation obtained after the Weierstrass substitution.

\section{Integration of \(I'(\lambda)\) and use of Coxeter's integrals}

We have
\[
I'(\lambda)=
\int_0^{\pi/2}
\frac{\cos^2\theta}
{(1+\lambda\cos\theta)\sqrt{(1+\lambda\cos\theta)^2-\cos^2\theta}}
\,d\theta.
\]
Integrating from 0 to \(\lambda>1\) yields the integral
\[
I(\lambda) - I(0) = \int_0^\lambda I'(s)\, ds
= \int_0^\lambda \int_0^{\pi/2}
\frac{\cos^2\theta}
{(1+s\cos\theta)\sqrt{(1+s\cos\theta)^2-\cos^2\theta}}
\,d\theta\, ds,
\]
which can be expressed, for $0<\lambda<2$ and using the standard reduction of quartic integrals with two real roots, in terms of incomplete elliptic integrals of the first kind $F$ and of the third kind $\Pi$ \cite{Byrd1971}:
\[
I'(\lambda)=\frac{2 i}{\sqrt{2-\lambda}\,\lambda(\lambda^2-1)}
\sqrt{\frac{1}{2+\lambda}}
\Bigg[
\lambda(1+\lambda)\,
F\!\left(
i\,\operatorname{arcsinh}\!\left(\frac{\sqrt{2-\lambda}}{\sqrt{\lambda}}\right)
\,\middle|\,
\frac{\lambda^2}{\lambda^2-4}
\right)
\]
\[
-\,2\Bigg(
(\lambda^2-1)\,
\Pi\!\left(
-\frac{\lambda}{\lambda-2};
\,i\,\operatorname{arcsinh}\!\left(\frac{\sqrt{2-\lambda}}{\sqrt{\lambda}}\right)
\,\middle|\,
\frac{\lambda^2}{\lambda^2-4}
\right)
\]
\[
+\,\Pi\!\left(
\frac{(\lambda-1)\lambda}{(\lambda-2)(1+\lambda)};
\,i\,\operatorname{arcsinh}\!\left(\frac{\sqrt{2-\lambda}}{\sqrt{\lambda}}\right)
\,\middle|\,
\frac{\lambda^2}{\lambda^2-4}
\right)
\Bigg)
\Bigg].
\]
More precisely, one has, for $-\frac{\pi}{2} < \phi < \frac{\pi}{2}$,
\[
F(\phi \mid m)
=
\int_{0}^{\phi}
\frac{d\theta}{\sqrt{1 - m \sin^2\theta}},
\]
and
\[
\Pi(n;\phi \mid m)
=
\int_{0}^{\phi}
\frac{d\theta}
{(1 - n \sin^2\theta)\,\sqrt{1 - m \sin^2\theta}}.
\]
For $\lambda=2$, one gets
\[
I(2) - I(0) = \int_0^2 I'(s)\, ds
= \int_0^2\int_0^{\pi/2}
\frac{\cos^2\theta}
{(1+s\cos\theta)\sqrt{(1+s\cos\theta)^2-\cos^2\theta}}
\,d\theta\, ds.
\]
Note that since the integrand is continuous and non-negative on compact subsets of $(0,2) \times [0,\pi/2]$, Fubini's theorem allows us to interchange the order of integration. It turns out that
$$
I(2)=A,
$$
the first Coxeter integral.

Although the integrands are different, Coxeter's evaluation shows that
$$
B=\pi^2/8=\int_0^{\pi/2}\arccos\!\left(\cos\theta\right)\,d\theta=\int_0^{\pi/2}\theta\,d\theta=I(0).
$$
Thus both Coxeter integrals $A$ and $B$ are encoded in the same family $I(\lambda)$: $A$ corresponds to $\lambda=2$, while $B$ appears at the boundary value $\lambda=0$. We thus have 
$$
\int_0^2 I'(s)\, ds
= \int_0^2 \int_0^{\pi/2}
\frac{\cos^2\theta}
{(1+s\cos\theta)\sqrt{(1+s\cos\theta)^2-\cos^2\theta}}
\,d\theta\, ds=A-B=\frac{5\pi^2}{24}-\frac{\pi^2}{8}
$$
yielding
$$
\int_0^2 \int_0^{\pi/2}
\frac{\cos^2\theta}
{(1+s\cos\theta)\sqrt{(1+s\cos\theta)^2-\cos^2\theta}}
\,d\theta\, ds=A-B=\frac{\pi^2}{12}.
$$

The third Coxeter integral $C$ can also be related to the same parametric framework after a simple transformation.
Recall that
\[
C=\int_0^{\pi/2}
\arccos\!\left(\frac{1-\cos\theta}{2\cos\theta}\right)d\theta.
\]
Writing
\[
\frac{1-\cos\theta}{2\cos\theta}
=
\frac{1}{2}\left(\sec\theta-1\right),
\]
and performing the change of variable $\theta\mapsto \frac{\pi}{2}-\theta$,
one obtains an integral involving
\[
\arccos\!\left(\frac{\sin\theta}{1+2\sin\theta}\right),
\]
which is structurally analogous to the integrand defining $I(\lambda)$, with $\cos\theta$ replaced by $\sin\theta$.
Thus $C$ may be viewed as another member of the same deformation scheme, corresponding to a rotated trigonometric parametrization. As in the cases of $A$ and $B$, Coxeter's evaluation yields
\[
C=\frac{11\pi^2}{72}.
\]

\section{Conclusion}

In this work, we have shown that classical Coxeter integrals, beyond their well-known exact evaluation as rational multiples of \(\pi^2\), reveal an elliptic structure underlying Coxeter's integrals. By embedding the first Coxeter integral $A$ into the one-parameter family 
\[
I(\lambda)=\int_{0}^{\pi/2}
\arccos\!\left(\frac{\cos\theta}{1+\lambda\cos\theta}\right)\,d\theta,
\]
and differentiating with respect to \(\lambda\), we obtained the explicit relation
\[
I'(\lambda)=
\int_0^{\pi/2}
\frac{\cos^2\theta}
{(1+\lambda\cos\theta)\sqrt{(1+\lambda\cos\theta)^2-\cos^2\theta}}
\,d\theta,
\]
or equivalently
\[
I'(\lambda) =
2 \int_{0}^{1} \frac{(1-t^2)^2}{(1+t^2)\left((1+\lambda)+(1-\lambda)t^2\right)\sqrt{(\lambda^2+2\lambda)+(4-2\lambda^2)t^2+(\lambda^2-2\lambda)t^4}} dt.
\]
We have also the expression
\begin{align}\label{ell}
I'(\lambda)&=\frac{2 i}{\sqrt{2-\lambda}\,\lambda(\lambda^2-1)}
\sqrt{\frac{1}{2+\lambda}}
\Bigg[
\lambda(1+\lambda)\,
F\!\left(
i\,\operatorname{arcsinh}\!\left(\frac{\sqrt{2-\lambda}}{\sqrt{\lambda}}\right)
\,\middle|\,
\frac{\lambda^2}{\lambda^2-4}
\right)\nonumber\\
&-\,2\Bigg(
(\lambda^2-1)\,
\Pi\!\left(
-\frac{\lambda}{\lambda-2};
\,i\,\operatorname{arcsinh}\!\left(\frac{\sqrt{2-\lambda}}{\sqrt{\lambda}}\right)
\,\middle|\,
\frac{\lambda^2}{\lambda^2-4}
\right)\nonumber\\
&+\,\Pi\!\left(
\frac{(\lambda-1)\lambda}{(\lambda-2)(1+\lambda)};
\,i\,\operatorname{arcsinh}\!\left(\frac{\sqrt{2-\lambda}}{\sqrt{\lambda}}\right)
\,\middle|\,
\frac{\lambda^2}{\lambda^2-4}
\right)
\Bigg)
\Bigg],
\end{align}
where $F$ is an incomplete elliptic integral of the first kind and $\Pi$ an incomplete elliptic integral of the third kind.

Integration then produces the following identities:
$$
\int_0^2 I'(\lambda)\,d\lambda=A-B=\frac{\pi^2}{12},
$$
with any of the three expressions of $I'(\lambda)$ mentioned above. Since the explicit elliptic representation (see Eq. (\ref{ell})) is valid only for
\[
0<\lambda<2,
\]
the integral over $[0,2]$ must be understood in the improper sense:
\[
\int_0^2 I'(\lambda)\,d\lambda
=
\lim_{\varepsilon\to0^+}
\int_{\varepsilon}^{2-\varepsilon} I'(\lambda)\,d\lambda.
\]
Two facts ensure that the limit exists. There is no genuine singularity in the original integral defining $I'(\lambda)$. The elliptic parametrisation degenerates at $\lambda=0$, but
$I'(\lambda)$ itself remains finite. Hence the limit as $\lambda\to0^+$ exists and is finite. The explicit formula produces a factor behaving like
\[
(2-\lambda)^{-1/2}.
\]
Thus, near $\lambda=2$,
\[
I'(\lambda)\sim \frac{D}{\sqrt{2-\lambda}}
\]
for some constant $D$. Since
\[
\int^{2}\frac{d\lambda}{\sqrt{2-\lambda}}
\]
is finite, the singularity is integrable. Therefore the improper integral converges.

This approach highlights that Coxeter integrals can serve as a bridge to relate elementary trigonometric integrals to elliptic integrals in a unified analytic framework. It opens the way for further explorations: other parametric deformations or generalizations of Coxeter-type integrals may produce additional integral identities involving elliptic integrals providing new insights into the interplay between classical analysis and elliptic functions.


\begin{thebibliography}{9}

\bibitem{Ahmed2002}
Z. Ahmed, Definitely an integral, {\it Amer. Math. Monthly} {\bf 109} (2002), 670--671. 

\bibitem{Byrd1971}
P.~F. Byrd and M.~D. Friedman,
\textit{Handbook of Elliptic Integrals for Engineers and Scientists},
2nd ed., Springer-Verlag, New York, 1971.

\bibitem{Borwein2004}
J. Borwein, D. Bailey, and R. Girgensohn, {\it Ahmed's Integral Problem}, §1.6 in Experimentation in Mathematics: Computational Paths to Discovery. Wellesley, MA: A K Peters, pp. 17-20, 2004.

\bibitem{Coxeter1937}
H.~S.~M.~Coxeter,
\newblock Some integrals connected with trigonometry,
\newblock {\em Math. Gaz.}, \textbf{21} (1937), 55--57.

\bibitem{Coxeter1968}
H.~S.~M.~Coxeter, {\it The Functions of Schl\"afli and Lobatschefsky}, in: Twelve geometric essays, Carbondale, Southern Illinois Univ. Press, 1968.

\bibitem{Coxeter1973}
H.~S.~M.~Coxeter, {\it Regular Polytopes}, 2nd edition, Dover, New York, 1973.

\bibitem{Nahin2015}
P. J. Nahin, {\it Inside Interesting Integrals}, Springer, 2015.
 
\bibitem{Schlafli1858}
L.~Schl\"afli,
\newblock Theorie der vielfachen Kontinuität,
\newblock in {\em Gesammelte mathematische Abhandlungen}, Vol.~I,
Birkhäuser, Basel, 1950 (original work 1858--1860).

\bibitem{Vidiani2003}
G.~Vidiani,
\newblock Les int\'egrales de Coxeter,
\newblock {\em Quadrature, Magazine de math\'ematiques pures et \'epic\'ees}, no.~50,
Oct.–Dec.~2003, pp.~7--12 (in French).

\end{thebibliography}
\end{document}